\input morse2.def

\document

\Kapitel\WhKapitel=6 The universal torsion class on the Whitehead space

In~\IBuch, Igusa defines the Whitehead space as the homotopy fibre
of a natural map from the stable homotopy of the classifying space~$BG$
of a subgroup~$G$ of the group of units~$R^\times$ of a ring~$R$
into the algebraic $K$-theory space~$BGL(R)^+$ of~$R$ itself.
In particular
	$$\Wh^h(R,G)\to Q(BG_+)\to \overline\Z\times BGL(R)^+$$
is a homotopy fibre sequence.
Here~$BG_+$ is the universal classifying space of~$G$
with an extra base point,
and~$Q$ denotes the stable homotopy functor~$\Omega^\infty S^\infty$.
The space~$\overline\Z\times BGL(R)^+\subset K(R)$
gives the algebraic $K$-theory of~$R$ in degrees~$\ge 1$,
where~$\overline\Z$ denote the image of~$\Z=K_0(\Z)$ in~$K_0(R)$,
$BGL(R)$ is the classifying space for the general linear
group~$GL(R)=\lim\limits_{\longrightarrow}GL(n,R)_\delta$,
and~$(\punkt)^+$ denotes Quillen's plus construction.
Note that all groups and rings here carry the discrete topology.

In Chapter~3 of~\IBuch,
Igusa constructs a model for the Whitehead space
out of ``systems of higher homotopies'',
which are precisely the ``combinatorial flat superconnections''
that we encountered in \CWBaseProp.
In Subsections~\NummerVon\SimpSubsection\ and~\NummerVon\SimpWhSubsection,
we review Igusa's construction of the simplicial Whitehead space
as the geometric realisation of a simplicial category~$\cWh^{\simp}(R,G)$.

Let~$p\colon M\to B$ be a proper submersion
equipped with a fibre-wise Morse function~$h\colon M\to\R$
and a flat vector bundle~$F\to M$.
We show in \FunctorSubsection\ that our construction
of a combinatorial flat superconnection in \CWSubsection\
gives a simplicial functor~$\Xi=\Xi_{p,h,F}$
from a certain simplicial category~$\Cal S$ associated to~$B$
to~$\cWh^{\simp}(R,G)$,
which is homotopy equivalent to the corresponding functor
constructed in~\IBuch, Chapter~4.

If the ring~$R$ is an $\R$-subalgebra of the ring~$\Cal M_n(\C)$
of $(n\times n)$-matrices over~$\C$,
and~$G\subset R^\times$ is any subgroup,
we develop two more models for the Whitehead space in \MixedWhSubsection,
given by differential superconnections and by mixed combinatorial
and differential superconnections.
We show that all these models are homotopy equivalent,
so that we may speak of ``the'' Whitehead space~$\Wh(R,G)$.
We want to point out that this equivalence of Whitehead categories
has been investigated independently by K. Igusa in~\ITwisted,
where he interprets both the map into the Whitehead space
and the bundle of the fibre-wise cohomology as $A_\infty$-functors.

In \ChernKapitel,
we will show that the form~$T_f(A',g^V,h^V)$ of \MVBTorsionDef\
is given by a universal cohomology class on a Whitehead space
if the family of cochain complexes~$(V,a_0)\to B$ is fibre-wise acyclic.
Since we may choose any of the models above for this Whitehead space,
this justifies calling the form~$T_f(A',g^V,h^V)$
a ``combinatorial torsion form'',
although it is defined in terms of a differential superconnection.
In \UniKapitel,
this construction will be generalised to the case that~$H=H^*(V,a_0)\to B$
carries a parallel metric.

\Abschnitt\SimpSubsection Simplicial sets and categories

We recall the definitions of simplicial sets and simplicial categories
as well as some elementary facts and examples needed later on.

\Definition\SimplSetDef
Let~$\Delta$ denote the small category
whose objects are the sets~$[k]=\{0,\dots,k\}$,
and whose morphisms are all mappings~$f\colon[j]\to[k]$
with~$f(i+1)\ge f(i)$ for all~$0\le i<j$.
If~$\Cal C$ denotes any category,
then a {\em simplicial object\/} in~$\Cal C$ is a contravariant
functor from~$\Delta$ to~$\Cal C$,
and a {\em simplicial morphism\/} between to simplicial objects in~$\Cal C$
is a natural transformation.

{\em Simplicial sets\/} and {\em simplicial maps\/}
are simplicial objects and morphisms in the category of sets and set maps.
{\em Simplicial categories\/} and {\em simplicial functors\/}
are simplicial objects and morphisms in the category
of small categories and covariant functors.
\enddefinition

If~$S$ is a simplicial set,
then~$S_k=S[k]$ is the {\em set of elements of degree}~$k$,
and~$Sf\colon S_k\to S_j$ a {\em (generalised) face operator\/}
for each morphism~$f\colon[j]\to[k]$ in~$\Delta$.
This terminology is compatible with \SimplexDef\ in~\TeilEins\
of a smooth simplicial complex,
where~$S_k$ was a set of smooth $k$-simplices~$\sigma\colon\Delta_k\to B$.
To comply with \SimplSetDef,
we must add degenerated simplices and degeneracy maps
to the complex defined in \SimplexDef.

If~$\Cal S$ is a simplicial category,
we still write~$\Cal S_k$ for the set of objects of degree~$k$.
The morphisms of~$\Cal S_k$ commute with
the face operators~$\Cal S_k\to\Cal S_l$.


The geometric realisation is a functor from the category
of simplicial sets to the category of topological spaces,
cf.~\Loday.
The geometric realisation of a small category is the realisation of its nerve,
which is again a simplicial set.
Similarly,
the geometric realisation of a simplicial category is the realisation
of its nerve,
which is now a bisimplicial set.

\Example\SimpCatExp
Let~$S$ be a smooth simplicial complex on~$B$
as in \SimplexDef.
Then we define a simplicial category~$\Cal S$:
\roster
\item
the set of objects of degree~$k$ is
	$$\Cal S_k
	=\bigl\{\,(\sigma,\sigma')
		\bigm|\sigma'\in S_k\text{ is a generalised face of }
			\sigma\in S_l\text{ for some }l\,\bigr\}\;;$$
\item
there is a unique morphism from~$(\sigma,\sigma')$
to~$(\tau,\tau')$
iff~$\sigma$ is a generalised face of~$\tau$
and~$\sigma'=\tau'$, 
and none otherwise;
\item
the generalised face operators~$f$ act by
	$$(\sigma,\sigma')\mapsto\(\sigma,\sigma'_{f(0),\dots,f(j)}\)\;.$$
\endroster
If~$S$ is sufficiently fine,
then the geometric realisation of~$\Cal S$
is homeomorphic to the set~$\overline B\subset B\times B$ of \BqSigmaDef.
In particular,
$B\cong\abs S$ is then a deformation retract of~$\abs{\Cal S}$.

Note that we have implicitly worked with~$\Cal S$ and~$\abs{\Cal S}$
in Sections~\NummerVon\GoodCoverSubsection\
and~\NummerVon\ComplexPatchSubsection.
We will describe the various complexes and superconnections
constructed in \ComplexKapitel\
as simplicial functors from~$\Cal S$
to several simplicial categories that are models
for Igusa's Whitehead space.
\endexample

\Abschnitt\SimpWhSubsection The simplicial Whitehead category

We start by recalling Igusa's simplicial Whitehead category~$\cWh^\simp(R,G)$
of~\IBuch, Chapter~3.
The typical object of this category is a
system of coefficients~$a_{\alpha\beta}(i_0,\dots,i_j)$ as in \CWBaseProp.
In this subsection,
$R$ is an associative ring with a unit,
and~$G\subset R^\times$ is any subgroup.

\Definition\GradedPosetDef
A {\em graded poset\/} ({\em graded p\/}artially {\em o\/}rdered {\em set\/})
is a finite $\Z$-graded set~$P=\bigcup_nP^n$
that is partially ordered by a relation~$\prec$
which is independent of the $\Z$-grading.
For~$\alpha\in P^n$, we set~$\ind(\alpha)=n$.
If~$P$, $Q$ are graded posets,
a {\em graded closed map from~$P$ to~$Q$\/}
is a map~$f\colon P\to Q$ such that for all~$\alpha$, $\beta\in P$,
$\ind\(f(\alpha)\)=\ind(\alpha)$,
and~$f(\alpha)\mathrel{\prec_Q}f(\beta)$ implies~$\alpha\prec_P\beta$.
\enddefinition

In particular,
a graded bijection~$P$ to~$Q$ is closed
iff the partial relation~$\prec_P$ is finer than the pullback of~$\prec_Q$.

\Definition\VectMatDef
Let~$P$ be a graded poset,
and let~$V$ be a 
left $R$-module.
Let~$V^P=\bigoplus_{\alpha\in P}V_\alpha$ be a direct sum of 
left $R$-modules with~$V_\alpha\cong V$ for all~$\alpha\in P$.
A {\em system of $V$-valued $P$-vectors~$v$ in degree}~$k$
assigns to each sequence~$0\le i_0\le\dots\le i_l\le k$
an element~$v(i_0,\dots,i_l)\in V$.
A {\em system of $P$-matrices~$a$ in degree}~$k$
assigns to each sequence~$0\le i_0\le\dots\le i_l\le k$
a matrix~$a(i_0,\dots,i_l)$ with coefficients in~$R$,
with rows and columns indexed by~$P$.
We call a system of matrices ({\em lower\/}) {\em triangular\/}
iff~$a_{\alpha\beta}(i_0,\dots,i_l)=0$ unless~$\beta\prec\alpha$.
Let~$V(k)$ and~$\Cal M_P\(R;k\)$
be the spaces of systems of $P$-vectors
and of $P$-matrices in degree~$k$,
and let~$T_P(R;k)\subset \Cal M_P(R;k)$
be the subspace of triangular systems.

We define a simplicial coboundary operator~$\delta$ on~$V(k)$ 
by
$$
	(\delta v)(i_0,\dots,i_l)
	=\sum_{j=0}^l(-1)^j\,v(i_0,\dots,\widehat{i_j},\dots,i_l)\;,
\tag1$$
and similarly on~$\Cal M_P(R;k)$.
We let~$\Cal M_P(R;k)$ act on~$V(k)$ 
by multiplication
$$
	(a*v)_\alpha(i_0,\dots,i_l)
	=\sum_{\beta\in P}\sum_{j=0}^l(-1)^{(l-j)\,(\ind\alpha-\ind\beta)}
		a_{\alpha\beta}(i_0,\dots,i_j)\,v_\beta(i_j,\dots,i_l)\;,
\tag2$$
and similarly,
we define a multliplication on~$\Cal M_P(R;k)$.
\enddefinition

\Remark\MatVectRem
The sign factor indicates that we have a ``horizontal'' grading,
given by the degree of the simplices we evaluate on,
and a ``vertical'' grading,
given by the grading of~$P$.
This is compatible with our orientation convention in \CWSubsection.
We leave it to the reader to check
that~``$*$'' is associative and that~``$\delta$''
is an odd derivation with respect to the total grading:
$$\alignat2
	(a* b)* v
	&=a*(b* v)\;,
	&\delta(a* v)
	&=(\delta a)* v+(-1)^{\deg a}\,a*(\delta v)\;,\\
	(a* b)* c
	&=a*(b* c)\;,
		\Mand
	&\delta(a* b)
	&=(\delta a)* b+(-1)^{\deg a}\,a*(\delta b)\;,\\
\endalignat$$
where we set
	$$\((-1)^{\deg a}a\)_{\alpha\beta}(i_0,\dots,i_j)
	=(-1)^{j+\ind\alpha-\ind\beta}\,a_{\alpha\beta}(i_0,\dots,i_j)\;.$$
\endremark

\Definition\SimpSupConnDef
If~$a$ is an element of~$\Cal M_P(R;l)$,
we call~$\delta+a$
a {\em simplicial $P$-super\-con\-nec\-tion of degree~$l$.\/}
We say that~$\delta+a$ is flat if
	$$(\delta+a*\mathord{})^2=0\tag1$$
as an operator on~$\Cal M_P(R;l)$.
We say that~$\delta+a$ is {\em of total degree~$1$}
if
	$$a_{\alpha\beta}(i_0,\dots,i_j)=0
		\Mtext{unless}
	\ind\alpha+j=\ind\beta+1\;.$$
We say that~$\delta+a$ is {\em triangular\/} if~$a$ is lower triangular
with respect to~$\prec_P$.

Let~$\delta+a$ be a flat simplicial $P$-superconnection.
If~$V=\bigoplus_{\alpha\in P}V_\alpha$ is a direct sum of left $R$-modules
and~$(C^*,\partial)$ is a cochain complex,
then a
{\em system of cochain maps~$I\colon(C,\partial)\to\(V(l),\delta+a\)$\/}
is an element~$I\in\(\Hom_\Z(C,V)\)^P(l)$,
such that
	$$(\delta+a*\mathord{})\circ I=I\circ\partial\;.\tag2$$
We say that~$I$ is {\em of degree~$k$,\/}
if~$I_\alpha(i_0,\dots,i_j)|_{C^n}=0$ unless~$\ind\alpha+j=n+k$.
\enddefinition

\Example\SimpSupConnExp
Let~$S$ be a smooth simplicial complex on~$B$,
and let~$a_{\alpha\beta}(\sigma)$ be defined as in \CWSubsection.
For each simplex~$\sigma\in S_l$,
let~$P$ be the set of leaves of~$\hat p\colon C(\abs\sigma)\to\abs\sigma$,
which is partially ordered by~$\prec_\sigma$.
If we set
	$$a_\sigma(i_0,\dots, i_k)
	=\cases
		a(\sigma_{i_0,\dots,i_k})
			&\text{if~$0\le i_0<\dots<i_k\le l$, and}\\
		0	&\text{otherwise,}
	\endcases
	\Formel\ExampleADef$$
then~$a\in T_P(\Z;l)$,
and~$\delta+a$ is a flat triangular
simplicial superconnection~$\delta+a_\sigma$
of total degree~$1$ by \CWBaseProp~(1) and~(3).

Assume that the bundle~$T^uX\to C$ is oriented,
and define~$I_\sigma\colon\Omega^*(M)\to\R^P(l)$ by
	$$I_{\sigma,\alpha}(i_0,\dots,i_j)\,\omega
	=\cases
		\int_{M_\alpha(\sigma_{i_0,\dots,i_j})}\omega
			&\text{if~$0\le i_0<\dots<i_k\le l$, and}\\
		0	&\text{otherwise.}
	\endcases
	\Formel\ExampleIDef$$
Then~$I_\sigma\colon\(\Omega^*(M),d^M\)\to\(\R^P(l),\delta+a\)$
is a system of cochain maps of degree~$0$
by Stokes' theorem and \CWBaseProp~(2).
\endexample

Let~$G$ be a subgroup of the multiplicative group~$R^\times$.
As a model for the Whitehead space mentioned at the beginning of this section,
we construct a simplicial category~$\cWh^\simp(R,G)$.
In fact,
we first define the pre-Whitehead category~$\cWh^{\prime\simp}(R,G)$.
The Whitehead category itself is obtained from~$\cWh^{\prime\simp}(R,G)$
by adding the so-called stabilisation morphisms,
and the (pre-\penalty10000)\hskip0pt\penalty0Whitehead
space is just the geometric realisation
of the (pre-)Whitehead category.
This category
is defined precisely as in \IBuch, Chapter~3.
However, we use a different,
but equivalent sign convention,
cf.~\IgusaSignRem.

\Definition\WhSimpDef
The {\em simplicial pre-Whitehead category~$\cWh^{\prime\simp}(R,G)$\/}
is defined as follows.
\roster
\item
The objects~$\cWh^{\prime\simp}_k(R,G)$ of degree~$k$ are pairs~$(P,a)$
where~$P$ is a graded poset,
and~$\delta+a$ is a flat triangular simplicial superconnection
of total degree~$1$ with~$a\in T_P(R;k)$,
such that~$a(i_0,\dots,i_l)=0$ unless~$i_0<\dots<i_l$.
\item
For~$(P,a)$,  $(Q,b)\in\cWh^{\prime\simp}_k(R,G)$,
the morphisms from~$(P,a)$ to~$(Q,b)$
are pairs~$(f,g)$
of a graded closed bijection~$f\colon P\rightarrow Q$
and a map~$g\colon P\rightarrow G$,
such that
	$$b_{f(\alpha),f(\beta)}\(i_0,\dots,i_l\)
	=g_\alpha\,a_{\alpha\beta}(i_0,\dots,i_l)\,g_\beta^{-1}$$
for all~$\alpha$, $\beta\in P$ and all~$0\le i_0<\dots<i_l\le k$.
\item
face operators are given by re-indexing:
	$$F_{j_0\dots j_m}(P,a)=(P,b)\in\cWh^{\prime\simp}_m(R;G)
		\Mtext{with}
	b(i_0,\dots,i_l)=a\(j_{i_0},\dots,j_{i_m}\)$$
for all sequences~$0\le j_0\le\dots\le j_m\le k$
and~$0\le i_0\le\dots\le i_l\le m$.
\endroster
The {\em simplicial pre-Whitehead
space\/}~$\Wh^{\prime\simp}(R,G)$ 
is the geometric realisation of the simplicial pre-Whitehead category.
\enddefinition

In order to keep the category~$\cWh^{\prime\simp}(R,G)$ small,
we agree that all our posets~$P$ are subsets of some large fixed set.
Since the poset~$P$ is mainly a book-keeping device,
this convention will not affect the homotopy type of any space
that we construct.

\Remark\WhFlatBdlRem
Fix a complex representation~$\rho\colon R\to \Cal M_r(\C)$.
Then the simplicial pre-Whitehead space carries two natural $\Z$-graded
flat complex vector bundles.
\roster
\item
There is a natural $\Z$-graded flat bundle~$\Cal V\to\Wh^{\prime\simp}(R,G)$
with fibre~$(\C^r)^P$ over each point of~$\abs{(P,a)}$
for each~$k$ and each~$(P,a)\in\cWh^{\prime\simp}_k(R,G)$.
The structure group of~$\Cal V$ is the group of $G$-monomial matrices.
\item
Because~$a(i)^2=0$,
we also have cohomology groups~$H^*\(V^P,a(i)\)$
for each~$0\le i\le k$.
Flatness of~$\delta+a$ implies
that~$\id+a(i_0,i_1)$ induces a cochain isomorphism from~$(V^P,a(i_1))$
to~$(V^P,a(i_0))$,
so we get isomorphisms
	$$\(\id+a(i_0,i_1)\)_*\colon H^*\(V^P,a(i_1)\)\to H^*\(V^P,a(i_0)\)$$
for~$0\le i_0\le i_1\le l$.
For~$i_0<i_1<i_2$,
the maps~$a(i_0,i_1,i_2)$ are cochain homotopies
between~$\(\id+a(i_0,i_1)\)\circ\(\id+a(i_1,i_2)\)$ and~$\(\id+a(i_0,i_2)\)$,
so all spaces~$H^*\(V^P,a(i)\)$ can be identified in a compatible way.
Morphisms~$(f,g)\colon(P,a)\to(Q,b)$ induce maps
	$$(f,g)_*(i)\colon H^*\(V^P,a(i)\)\to H^*\(V^Q,b(i)\)$$
that are compatible with the isomorphisms~$\(\id+a(i_0,i_1)\)_*$ above.
This way,
we obtain a natural $\Z$-graded
flat bundle
	$$\Cal H
	=H^*(\Cal V,a)
	\to\Wh^{\prime\simp}(R,G)$$
of cohomology groups.
Note that the structure group of~$\Cal H$
is a subgroup of~$GL(R)$,
which is in general not related to~$G$.
\endroster
\endremark

So far,
each connected component of~$\Wh^{\prime\simp}(R,G)$ is built
from graded posets~$P$ having the same number~$n_k$ of elements degree~$k$
for each~$k$.
Using Igusa's stabilisation construction,
we add more morphisms to the underlying categories~$\cWh^{\prime\simp}(R,G)$
between connected components of~$\Wh^{\prime\simp}(R,G)$.
This gives us the (stabilised) Whitehead space~$\Wh^\simp(R,G)$.
These extra morphisms are necessary to account for the changes
of the family Thom-Smale complex
that occur at birth-death singularities.
We will check that all our constructions are compatible with stabilisation.
In particular,
we prove that the universal torsion class~$T(R,G)$
we are going to define in \ChernAltSubsection\ extends to~$\Wh(R,G)$.

The only new ingredients in the final construction of~$\Wh(R,G)$
are elementary expansions of simplicial superconnections.
Let~$P\dg j$ denote the poset consisting of two elements~$\alpha^-$ in
degree~$j$ and~$\alpha^+$ in degree~$j+1$,
with~$\alpha^-\prec\alpha^+$.
Let~$a\dg j\in \Cal M_{P\dg j}(R,k)$ denote the system of $P\dg j$-matrices
given by
	$$a\dg j(i_0,\dots,i_l)
	=\cases
		\bigl(\smallmatrix 0&1\\&0\endsmallmatrix\bigr)
			&\text{if~$l=0$, and}\\
		\bigl(\smallmatrix 0\\&0\endsmallmatrix\bigr)
			&\text{if~$l\ge1$.}
	\endcases\Formel\ElExADef$$
Then~$\bigl(P\dg j,a\dg j\bigr)$ is an object in~$\cWh_k^{\prime\simp}(R,G)$,
which we call the {\em elementary object of vertical degree~$j$
in~$\cWh_k^{\prime\simp}(R,G)$.\/}

We can define a direct sum on the objects of~$\cWh_k^{\prime\simp}(R,G)$
by setting
	$$(P,a)\oplus(P',a')=(P\cup P',a\oplus a')\;.$$
Here, the partial ordering on~$P\cup P'$
is the union of the partial orderings on~$P$ and~$P'$,
so~$\alpha\prec\beta$
iff either~$\alpha$, $\beta\in P$ and~$\alpha\prec\beta$ in~$P$,
or~$\alpha$, $\beta\in P'$ and~$\alpha\prec\beta$ in~$P'$.

\Definition\ElExDef
Let~$(P,a)$ be an object of~$\cWh_k^{\prime\simp}(R,G)$.
Then the object~$(P,a)\oplus\bigl(P\dg j,a\dg j\bigr)$
is called the {\em elementary expansion of\/~$(P,a)$ in degree\/~$j$.\/}
\enddefinition

We can now define the stabilised simplicial Whitehead category.

\Definition\WhStabDef
The {\em (stabilised) simplicial Whitehead category\/}~$\cWh^{\simp}(R,G)$
has the same objects and face operators as~$\cWh^{\prime\simp}(R,G)$.
Its morphisms are iterated formal compositions
of~$\cWh^{\prime\simp}(R,G)$-morphisms and elementary expansions
	$$(P,a)\to(P,a)\oplus\(P\dg j,a\dg j\)\;.$$
The {\em simplicial Whitehead space\/}~$\Wh^{\simp}(R,G)$
is the geometric realisation of the
simplicial Whitehead category.
\enddefinition


\Remark\WhFlatBdlRemTwo
Note that there is a natural inclusion
	$$\Wh^{\prime\simp}(R,G)\to\Wh^{\simp}(R,G)\;,$$
since~$\cWh^\simp(R,G)$ has more morphisms than~$\cWh^{\prime\simp}(R,G)$,
and thus~$\Wh^{\simp}(R,G)$ has additional bi-simplices.
We will call a flat bundle or a cohomology class on~$\cWh^{\prime\simp}(R,G)$
{\em compatible with stabilisation\/}
if it is the pull-back of a corresponding object on~$\Wh^{\simp}(R,G)$.

The morphisms of the stabilised category still induce
isomorphisms on the level of cohomology,
because expansion pairs do not contribute to the cohomology.
In particular,
we still have a well-defined bundle~$\Cal H\to\Wh^\simp(R,G)$
as in \WhFlatBdlRem~(2),
which pulls back to~$\Cal H$.
However,
the vector bundle~$\Cal V\to\Wh^{\prime\simp}(R,G)$
extends to~$\Wh^{\simp}(R,G)$ only as an equivalence class
of vector bundles up to elementary expansions.
This means that we will not define higher torsion classes directly
on~$\Wh^{\simp}(R,G)$,
but we will instead look for classes on~$\Cal H\to\Wh^{\prime\simp}(R,G)$
that are compatible with stabilisation.
\endremark

We also recall the definition of the acyclic Whitehead space,
which plays an important role in~\IBuch.

\Definition\WhAcyclicDef
The {\em acyclic simplicial Whitehead category\/}~$\cWh^{\simp,h}(R,G)$
is the full subcategory of~$\cWh^{\simp}(R,G)$
containing all objects~$(P,a)$ such that the complex~$(V^P,a(0))$ is acyclic.
The {\em acyclic simplicial Whitehead space\/}~$\Wh^{\simp,h}(R,G)
=\bigl|N\cWh^{\simp,h}(R,G)\bigr|$ is the geometric realisation
of this category.
\enddefinition

We remark that~$\Wh^{\simp,h}(R,G)$
is the union of all connected components of~$\Wh^{\simp}(R,G)$
for which the bundle~$\Cal H$ vanishes.

\Abschnitt\FunctorSubsection A functor to the Whitehead category

Let~$M\to B$ be a bundle of compact manifolds,
and let~$F\to M$ be a flat bundle
with structure group~$G\subset R^\times$,
such that~$\rho|_G\colon G\to GL_r(\C)$ is faithful (injective).
E.g., if~$F$ is a complex vector bundle,
we may take~$R=\Cal M_r(\C)$ to be the ring of complex $r\times r$-matrices,
and~$G\subset GL_rR$,
where~$r$ is the complex rank of~$F$.
If~$F$ is unitarily flat, then~$G\subset U(r)$.
In general,
if~$F$ is associated to an $r$-dimensional complex representation
of~$\pi_1(M)$ that factors over a group~$G$,
so
	$$\pi_1(M)@>>>G@>\rho'>>GL_r(\C)\;,$$
we may take~$R=\Z(G)$ the group ring of~$G$,
and the representation~$\rho\colon\Z(G)\to \Cal M_r(\C)$ induced by~$\rho'$.

Let~$h$ be a fibre-wise Morse function on~$M$ with fibre-wise critical set~$C$,
and let~$S$ be a smooth simplicial structure on~$B$
as in \SimplexDef, which is fine enough to satisfy~\SimplFineFormel.
We want to interpret the construction in Sections~\NummerVon\CWSubsection\
and~\NummerVon\GoodCoverSubsection\ of~\TeilEins\
as the construction of a functor from the simplicial category~$\Cal S$
of \SimpCatExp\ to~$\cWh^{\prime\simp}(R,G)$.
Throughout this subsection,
we fix a ring~$R$ and a complex representation~$\rho\colon R\to \Cal M_r(\C)$.
Then there are natural flat complex vector bundles~$\Cal V$,
$\Cal H\to\Wh^{\prime\simp}(R,G)$ as in \WhFlatBdlRem.
Note that unless~$T^uX\to C$ is oriented,
we have to make sure that~$-1\in\rho'(G)$ by enlarging~$G$ if necessary.

Supposing that~$M$ is connected,
let us fix a basepoint~$o\in M$.
For each~$\sigma\in S$,
let~$P(\sigma)$ denote
the set of connected components of~$C|_{\inn M_\alpha(\sigma)}$,
which we order by the partial relation~$\prec_\sigma$ of~\CWPrecDef.
If~$\tau$ is a face of~$\sigma$,
then the relation~$\prec_\sigma$ is finer than~$\prec_\tau$,
so the natural identification
	$$f_{\tau,\sigma}\colon P(\sigma)\to P(\tau)$$
is a closed bijection.

We take the CW structure on~$M$
with cells~$M_\alpha(\sigma)$ for each simplex~$\sigma\in S$
and each fibre-wise critical level~$C_\alpha$ of~$h$ over~$\abs\sigma$
as in \CWSubsection.
For each~$\alpha\in P(\sigma)$,
we fix a path~$\gamma_\alpha(\sigma)$ from~$o$
to some point in~$\inn M_\alpha(\sigma)$.
Let us identify the space~$\Cal O_F\(M_\alpha(\sigma)\)$ of parallel sections
of~$F|_{\inn M_\alpha(\sigma)}$ with~$F_o$
by parallel translation along~$\gamma_\alpha(\sigma)$.
We also fix an orientation~$o_\alpha(\sigma)$ of~$T^uX|_{C_\alpha(\inn\sigma)}$
and set
	$$V(\sigma)=\bigoplus_\alpha V_\alpha(\sigma)\;.$$

Working with the simplicial set~$\Cal S$ of \SimpCatExp,
for each~$\sigma\in\Cal S_k$
we have to construct an element~$a\in T_{P(\sigma)}(R;k)$
such that~$\delta+a$ is a flat simplicial superconnection of total degree~$1$.
Therefore,
we fix orientations on the bundles~$T^uX|_{C_\alpha(\inn\sigma)}$
and regard the coboundary operator~$\delta+a$ constructed in \CWBaseProp\
for trivial coefficients,
and in \GCIntCor\ for the general case.
The coefficient~$a_{\alpha\beta}(\sigma_{i_0\dots i_j})$
arises as a sum of contributions of flow lines of~$W$
going from~$C_\alpha(\sigma_{i_0})$ to~$C_\beta(\sigma_{i_j})$.
Together with~$\gamma_\alpha(\sigma)$, $\gamma_\beta(\sigma)$
and paths inside~$C_\alpha(\abs\sigma)$ and~$C_\beta(\abs\sigma)$,
these connecting lines give rise to a loop in~$M$.
Since parallel translation along this loop acts by an element of~$G$,
we may regard~$a_{\alpha\beta}(\sigma_{i_0\dots i_j})$
as an element of~$\Z(G)$,
or of any ring~$R$ containing~$G$.
By \CWBaseProp~(1), (3) and by~\GCIntCor~(2), (3),
these coefficients give rise to a flat triangular
simplicial superconnection~$\delta+a_\sigma$ of total degree~$1$.

Now to each pair~$(\sigma,\sigma')\in\Cal S$,
where~$\sigma'=\sigma_{i_0\dots i_j}$ is a face (or degeneracy) of~$\sigma$,
we associate the graded poset~$P(\sigma)$
and the restriction of the simplicial superconnection~$\delta+a_\sigma$
to~$\sigma'$ following \WhSimpDef~(3).
Set
	$$\Xi^{\simp}_{p,h,F}(\sigma,\sigma')
	=\(P(\sigma),a_\sigma|_{\sigma'}\)
		\quad\in\cWh^{\prime\simp}_j(R,G)\;.
	\Formel\WhXiDef$$
Face operators in~$\Cal S$ map~$(\sigma,\sigma')$
to~$(\sigma,\sigma'')$,
where~$\sigma''$ is a face of~$\sigma'$.
This implies that~$\Xi^{\simp}_{p,h,F}$
commutes with generalised face operators.

\Proposition\WhFunctorProp.
The assigment~$\Xi^{\simp}_{p,h,F}$ constructed above
is a simplicial functor from~$\Cal S$
to the simplicial pre-Whitehead category~$\cWh^{\prime\simp}(R,G)$.
The flat bundles~$V\to B$ and~$H\to B$
are naturally isomorphic to the pull-backs by~$\abs{\Xi_{p,h,F}}$
of the bundles~$\Cal V$, $\Cal H\to\Wh^{\prime\simp}(R,G)$ of \WhFlatBdlRem.
If~$S'$ is another sufficiently fine smooth simplicial complex on~$B$,
$\Cal S'$  is the corresponding simplicial category,
and~$\Xi^{\simp\prime}_{p,h,F}\colon\Cal S'\to\cWh(R;G)$ is another functor
constructed as above,
possibly with different choices than before,
then the diagramme
$$\psset{linewidth=0.5pt,nodesep=5pt}
	\matrix
		&&\Rnode{WhFcS}{\abs{\Cal S}}\\\vphantom{\Big|}\\
		\Rnode{WhFB}{B}&&
		&&\Rnode{WhFcWh}{\hphantom{B}\llap{$\Cal W$}}h^{\prime\simp}(R;G)\\
		\vphantom{\Big|}\\&&\Rnode{WhFcSx}{\abs{\Cal S'}}
	\endmatrix
\ncline{->}{WhFB}{WhFcS}\aput[0pt]{:U}{\sim}
\ncline{<-}{WhFcSx}{WhFB}\aput[0pt]{:U}{\sim}
\ncline{->}{WhFcS}{WhFcWh}\Aput{\scriptstyle\Xi^\simp_{p,h,F}}
\ncline{->}{WhFcSx}{WhFcWh}\Bput{\scriptstyle\Xi^{\simp\prime}_{p,h,F}}$$
commutes up to homotopy,
where the maps on the left are constructed using partitions of unity
as in~\PartUnFormel.

\Remark\WhSimpUniversalRem
Conversely,
given~$S$ and~$\Cal S$ as before,
a functor~$\Xi\colon\Cal S\to\cWh^{\prime\simp}$
induces
\roster
\item
a flat $\Z$-graded complex vector bundle~$(V,\nabla^V)\to B$
of the form~$V=\hat p_*\widehat V$,
where~$\widehat V\to C$ is a vector bundle of rank~$r$
with structure group~$G$ on some finite cover~$\hat p\colon C\to B$ of~$B$;
and
\item
a flat simplicial superconnection~$\delta+a$ of total degree~$1$
acting on the space of $V$-valued simplicial cochains,
which is triangular over each simplex~$\sigma$ with respect to some partial
ordering on the leaves of~$C|_{\inn\sigma}$.
\endroster
In other words,
the coefficients~$a_{\alpha\beta}(\sigma)$ satisfy \CWBaseProp~(1) and~(3).
Note that property~(2) in \CWBaseProp\ is meaningless here because we cannot
reconstruct the space~$M$ from~$V$ and~$a$.

This means that~$\Wh^\simp(R,G)=|\cWh^\simp(R;G)|$ is a ``universal space''
for the flat triangular simplicial superconnections of total degree~$1$
described above.
The geometric realisation
	$$\xi^\simp_{p,h,F}=\bigl|\Xi^\simp_{p,h,F}\bigr|
		\colon B\to\Wh^\simp(R,G)$$
is well-defined up to homotopy;
it is the ``classifying map''
for the given flat triangular simplicial superconnection.
\endremark

\demo{Proof\/ \rm of \WhFunctorProp}
Recall that there is an $\Cal S$-morphism from~$(\sigma,\sigma')$
to~$(\tau,\tau')$ iff~$\tau$ is a face of~$\sigma$ and~$\sigma'=\tau'$,
and that this morphism is unique.
We noted above that the natural identification~$f_{\tau,\sigma}
\colon P(\sigma)\to P(\tau)$ is a closed bijection of graded posets.

For~$\alpha\in P(\sigma)$,
the trivialisations of~$V_\alpha(\sigma)$ and~$V_{f_{\tau,\sigma}\alpha}(\tau)$
differ by an element~$g$ of~$G$.
If the orientations on~$C_\alpha(\sigma)$ chosen above restricts
to the one chosen on~$C_{f_{\tau,\sigma}\alpha}(\sigma)$,
we set~$g_{\tau,\sigma;\alpha}=g$,
otherwise we set~$g_{\tau,\sigma;\alpha}=-g$.

Note that by \CWBaseProp~(1),
$a_{\sigma;\alpha,\beta}(\sigma'')=0$
unless~$f_{\tau,\sigma}\beta\prec_\tau f_{\tau,\sigma}\alpha$.
Now it is easy to check that
	$$a_{\tau;f_{\tau,\sigma}\alpha,f_{\tau,\sigma}\beta}(\sigma'')
	g_{\tau,\sigma;\alpha}\circ a_{\sigma;\alpha,\beta}(\sigma'')
	\circ g_{\tau,\sigma;\beta}^{-1}$$
for all generalised faces~$\sigma''$ of~$\sigma'$.
Thus, the pair
	$$\Xi^\simp_{p,h,F}\((\sigma,\sigma')\to(\tau,\sigma')\)
	=\(f_{\tau,\sigma},g_{\tau,\sigma}\)$$
is a morphism in~$\cWh^{\prime\simp}(R;G)$.
Functoriality of~$\Xi^\simp_{p,h,F}$ is now trivial.

To show that~$\xi^\simp_{p,h,F}$ is well-defined up to homotopy,
we regard two smooth simplicial complexes on~$B\times\{0\}$
and~$B\times\{1\}$.
These can be combined to a smooth simplicial complex
on~$\overline B=B\times[0,1]$
with corresponding simplicial category~$\overline{\Cal S}$.
We may now construct a functor
	$$\overline\Xi^\simp_{p,h,F}
	\colon\overline{\Cal S}\to\cWh^{\prime\simp}(R,G)$$
that restricts to given functors on simplices in~$B\times\{0\}$
and~$B\times\{1\}$.
%
%
Our claim follows from the relativity of the constructions
in \ComplexKapitel, see e.g.\ the proof of \RelativityProp.
\qed\enddemo

Next,
we want to show
that the functor~$\Xi^\simp_{p,h,F}\colon\Cal S\to\cWh^{\prime\simp}(R,G)$
constructed above
is related to the one constructed in~\IBuch, Section~4.3,
by a natural transformation.
This will then imply that the realisations are homotopic.

We therefore go back to \CWSubsection,
where we constructed a cell complex
consisting of one cell~$M_\alpha(\sigma)$ for each simplex~$\sigma$
of a simplicial decomposition~$S$ of~$B$ and each fibre-wise
critical level~$\alpha$ of~$h$ over~$\inn\sigma\subset B$.
We may assume that~$S$ is sufficiently fine for the constructions
both in \ComplexKapitel\ and in Chapter~4 of~\IBuch,
e.g., we may take~$S$ to be the first barycentric subdivision
of a simplicial complex satisfying~\CWFineAss\ and~\SimplFineFormel.
This implies that any sequence
of simplices~$\sigma\supset\tau_0\supset\dots\supset\tau_k$ in~\IBuch\
corresponds to a single simplex in \ComplexKapitel.

The next technical issue is the choice of partial ordering.
In~\ComplexKapitel, we work with the partial order~$\prec_\sigma$
of~\CWPrecDef,
while Igusa uses a ``standard'' and a ``minimal order.
The standard order is defined by saying that
	$$\alpha<^s_\sigma\beta\Mtext{iff}h_\alpha<h_\beta
		\Rtext{on~$\abs\sigma$.}$$
The minimal order is generated by~$\alpha<^{\min}\beta$
if there exists a negative gradient flow line from~$C_\beta(\sigma)$
to~$C_\alpha(\sigma)$ for a fixed vertical metric~$g\TX$.

Note that both~$\prec_\sigma$ and~$<^{\min}_\sigma$ are refinements
of the corresponding orders on the faces of~$\sigma$,
while~$<^s_\sigma$ behaves the other way.
Also, the standard order is always finer than~$\prec$,
which is in turn finer than the minimal order by \CWPrecRemark.
If~$(P,a)$ is an object of~$\cWh^{\simp}(R,G)$,
then there is a morphism from~$(P,a)$ to the object~$(P^{\min},a)$,
where~$P^{\min}$ is the same set as~$P$,
but equipped with the minimal ordering with respect to~$a$.
Let~$\Xi^{\min}_{p,h,F}\colon\Cal S\to\cWh^{\prime\simp}(R,G)$
be the corresponding functor,
then there is a natural transformation from~$\Xi^\simp_{p,h,F}$
to~$\Xi^{\min}_{p,h,F}$,
so their geometric realisations are homotopic.

But the functor~$\Xi^{\min}_{p,h,F}$ gives an admissible choice
of elements~$\Xi^{\min}(\sigma)\in \Cal M_{P^{\min}(\sigma)}(R;k)$
for the construction leading to Igusa's Theorem~4.3.8
(recall that~$\sigma$ is a simplex of the barycentric subdivision,
so it corresponds to a sequence of simplices in~\IBuch).
In particular,
the realisation of the
functor~$\Xi^{\simp}_{p,h,F}\colon\Cal S\to\cWh^{\prime\simp}(R,G)$
is homotopic to the pullback of the realisation
of Igusa's~$\Xi^{\min}$ to~$\abs{\Cal S}$ by Theorem~4.3.8 in~\IBuch.
%
Let us summarize this.

\Theorem\IgusaFunctorThm.
The constructions in Chapter~4 of~\IBuch\ and in~\ComplexKapitel\ of this paper
give rise to homotopic maps from~$\abs{\Cal S}\sim B$
to~$\Wh^{\prime\simp}(R,G)$.

\Abschnitt\MixedWhSubsection
The differential and the mixed model of the Whitehead space

We introduce two more models~$\cWh^\diff(R,G)$
and~$\cWh^\mix(R,G)$ of the Whitehead space
(note that our notion of~$\cWh^\diff(R,G)$ differs from Igusa's). 
Because we have to interpolate and differentiate the coefficients
of mixed and differential superconnections,
the ring~$R$ must now be an $\R$-algebra
(typically, $R=\Cal M_r(\Bbbk)$ for~$\Bbbk=\R$, $\C$ or~$\Bbb H$).
However,
we may still work with any suitable group~$G\subset R^\times$.
We still fix a complex representation~$\rho\colon R\to \Cal M_r(\C)$.

Let~$\Delta_k$ denote the standard affine $l$-simplex,
and let~$f_{i_0\dots i_l}\colon\Delta_l\to\Delta_k$ be the unique
affine map sending the $j$-th vertex of~$\Delta_l$ to the $i_j$-th vertex
of~$\Delta_k$.
Let~$T_P(R)\subset \Cal M_P(R)$ denote the algebra of $P$-triangular matrices,
cf.\ \VectMatDef.
For a superconnection~$A'$
on the trivial bundle~$V^P\times\Delta_k\to\Delta_k$,
we write
	$$A'=d+\sum_ja'_j$$
with~$a'_j\in\Omega^j(\Delta_k;\End V^P)$ as usual.
Recall that superconnections of total degree~$1$ were defined in \TotDegOneDef.

\Definition\WhDiffDef
The {\em differential pre-Whitehead category\/}~$\cWh^{\prime\diff}(R,G)$
is defined as follows.
\roster
\item The objects~$\cWh_k^{\prime\diff}(R,G)$ of degree~$k$ are pairs~$(P,a')$
where~$P$ is a graded poset,
and~$d+a'$ is a flat superconnection of total degree~$1$
on the bundle~$\Delta_k\times(\C^r)^P\to\Delta_k$
with~$a'\in\Omega^*\(\Delta_k,T_P(R)\)$.
\item
For~$(P,a')$ and~$(Q,b')\in\cWh_k^{\prime\diff}(R,G)$,
the morphisms from~$(P,a')$ to~$(Q,b')$ are pairs~$(f,g)$
of a graded closed bijection~$f\colon P\to Q$
and a map~$g\colon P\to G$,
such that
	$$b'_{f(\alpha),f(\beta)}
	=g_\alpha\,a'_{\alpha\beta}\,g_\beta^{-1}
		\in\Omega^*\(\Delta_k,R\)
	\Rtext{for all~$\alpha$, $\beta\in P$.}$$
\item
Generalised
face operators~$F_{i_0\dots i_l}$
act by pulling~$d+a'$ back
by the corresponding affine map~$f_{i_0\dots i_l}\colon\Delta_l\to\Delta_k$,
so
	$$F_{i_0\dots i_l}(P,a')
	=\(P,f_{i_0\dots i_l}^*a'\)\;.$$
\endroster
The {\em differential
pre-Whitehead space\/}~$\Wh^\diff(R,G)$
is the geometric realisation of the differential pre-Whitehead category.
%
\enddefinition

\Remark\WhDiffFlatBdlRem
As in \WhFlatBdlRem,
the complex representation~$\rho$ gives rise
to a flat bundle~$\Cal V\to\Wh^{\prime\diff}(R,G)$.
Over each geometric simplex of~$\Wh^{\prime\diff}(R,G)$,
we consider the space of smooth $\Cal V$-valued forms.
Define~$\Omega^*\(\Wh^{\prime\diff}(R,G),\Cal V\)$ to be the set of
functions that assign a smooth $\Cal V$-valued form~$\omega(\sigma)$
to each simplex~$\sigma$ of~$\Wh^{\prime\diff}(R,G)$,
such that~$\omega(\sigma')=\omega(\sigma)|_{\sigma'}$
for each face~$\sigma'$ of~$\sigma$.
Then~$d+a'$ induces a piece-wise smooth flat superconnection~$\Cal A'$
on~$\Omega^*\(\Wh^{\prime\diff}(R,G),\Cal V\)$.
Let~$\Cal H=H^*(\Cal V,a'_0)\to\Wh^{\prime\diff}(R,G)$ be the flat bundle
of the fibre-wise cohomology
of~$\(\Omega^*(\Wh^{\prime\diff};\Cal V),\Cal A'\)$,
then~$\Cal H$ also carries a natural flat connection,
which is induced by~$\Cal A'$.
\endremark

Arguing as in \FunctorSubsection,
the differential pre-Whitehead space
is universal for the kind of $S$-piece-wise smooth flat superconnections
on vector bundles over simplicial complexes
that we discussed at the beginning of \SmoothingSubsection.
Similarly,
we now define the mixed pre-Whitehead space,
which is universal for the mixed superconnections
constructed in \SimplPatchProp.
Let~$\Delta_{k\setminus j}$
denote the span of the vertices~$j$, $j+1$, \dots, $k$ of~$\Delta_l$.

\Definition\MixedMatrixDef
Let~$P$ be a graded poset.
A {\em mixed system of triangular $P$-matrix-valued forms\/} on~$\Delta_k$
is a pair~$(a,a')$,
where~$a\in T_P(R;k)$ is a triangular system of $P$-matrices,
and~$a'$ assigns a differential form~$a'(i_0,\dots,i_{l-1})
\in\Omega^*\(\Delta_{k\setminus i_{l-1}};T_P(R)\)$
to each sequence~$0\le i_0\le\dots\le i_{l-1}\le k$ with~$0\le l$,
with
	$$a'(i_0,\dots,i_{l-1})|_{\Delta_k(i_l)}=a(i_0,\dots,i_l)\tag1$$
if~$i_{l-1}\le i_l\le k$.
Let~$M^\mix_P(\Delta_k;R)$ denote the space of all
mixed systems of triangular $P$-matrix-valued forms with coefficients in~$R$.

We define two operators~$\delta$ and~$d$
on~$M^\mix_P(\Delta_k;R)$ by~$\delta(a,a')=(\delta a,\delta'a')$
and~$d(a,a')=(0,d'a')$ with
$$\align
	(\delta'a')(i_0,\dots,i_{l-1})
	&=\sum_{j=0}^{l-1}(-1)^j\,a'\(i_0,\dots,\widehat{i_j},\dots,i_{l-1}\)
		+(-1)^l\,a\(i_0,\dots,i_{l-1}\)\;,\tag2\\
	(d'a')(i_0,\dots,i_{l-1})
	&=(-1)^ld\(a'(i_0,\dots,i_{l-1})\)\;.\tag3
\endalign$$
We also define a multiplication on~$M^\mix_P(\Delta_k;R)$
by~$(a,a')*(b,b')=(a*b,a'*'b')$,
with
$$\multline
	(a'*'b')_{\alpha\gamma}(i_0,\dots,i_{l-1})\\
	=\sum_{j=0}^{l-1}\sum_{\beta\in P}(-1)^{(l-j)(\ind\alpha-\ind\beta)}
			\,a_{\alpha\beta}(i_0,\dots,i_j)
			\,b'_{\beta\gamma}(i_j,\dots,i_{l-1})\\
		+\(a'(i_0,\dots,i_{l-1})
			\,b'(\emptyset)\)_{\alpha\gamma}\;.
\endmultline\tag4$$
In~(2) and~(4),
$a\(i_0,\dots,i_j\)$ is regarded
as a constant $T_P(R)$-valued function on~$\Delta_{k\setminus i_j}$.
\enddefinition

\Remark\MixedMatrixRem
We collect some properties of the definition above.
\roster
\item
Note that~(2), (3) and~(4) are compatible with~(1).
Note also that by~(1),
it would have been enough to specify~$a'$ instead of~$(a,a')$,
this is why we may write~$\delta'a'$ instead of~$\delta'(a,a')$
and~$a'*'b'$ instead of~$(a,a')*'(b,b')$.
We have chosen the notation~$(a,a')$ mainly for compatibility
with \ComplexPatchSubsection,
see also \SimplPatchRem.
\item
We remark in particular
that~$a'(\emptyset)\in\Omega^*\(\Delta_k;\End V^P\)$,
so that we can define a superconnection~$A'=d+a'(\emptyset)$
on the trivial bundle with fibre~$V^P$.
Let~$a'_j=a'(\emptyset)\dg j\in\Omega^j\(\Delta_k;\End V^P\)$
denote its component of degree~$j$.
\item
The sign factors indicate that we have again two gradings:
first,
the simplicial degree,
which is~$l$ for~$a'(i_0,\dots,i_{l-1})$
and~$a(i_0,\dots,i_l)$;
second,
we have the sum of the exterior degree and the degree given by~$P$
on~$\Omega^*\(\Delta_{k\setminus i_{l-1}};T_P(R)\)$.
Thus, $(a,a')\in \Cal M_P^\mix(\Delta_k;R)$ is of total degree~$m$ iff
	$$a'_{\alpha\beta}(i_0,\dots,i_{l-1})
	\in\Omega^{m-l-\ind\alpha+\ind\beta}
		\(\Delta_{k\setminus i_{l-1}}\)\;.$$
Note that~``$d$'' and~``$\delta$'' are super-commuting odd derivations
with respect to the total degree,
and that~``$*$'' is associative.
\endroster
\endremark

\Definition\WhMixedDef
The {\em mixed pre-Whitehead category\/}~$\cWh^{\prime\mix}(R,G)$
is defined as follows.
\roster
\item The objects~$\cWh_k^{\prime\simp}(R,G)$ of degree~$k$
are triples~$(P,a,a')$,
where~$P$ is a graded poset,
and~$(a,a')\in \Cal M_P^\mix(\Delta_k;R)$ is of total degree~$1$
and satisfies
	$$a(i_0,\dots,i_l)=a'(i_0,\dots,i_{l-1})=0$$
unless~$i_0<\dots<i_l$,
and
	$$d(a,a')+\delta(a,a')+(a,a')*(a,a')=0\;.$$
\item
For~$(P,a)$ and~$(Q,a')\in\cWh_k^{\prime\mix}(R,G)$,
the morphisms from~$(P,a)$ to~$(Q,a')$ are pairs~$(f,g)$
of a graded closed bijection~$f\colon P\to Q$
and a map~$g\colon P\to G$,
such that
	$$a'_{f(\alpha),f(\beta)}(i_0,\dots,i_{l-1})
	=g_\alpha\,a_{\alpha\beta}(i_0,\dots,i_{l-1})\,g_\beta^{-1}
		\qquad\in\Omega^*\(\Delta_{k\setminus i_{l-1}},R\)$$
for all~$\alpha$, $\beta\in P$
and all~$0\le i_0\le\dots\le i_{l-1}\le k$.
\item
For a sequence~$0\le j_0\le\dots\le j_m\le k$,
the face operator~$F_{j_0\dots j_m}
\colon\cWh_k^{\prime\mix}(R,G)\to\cWh_m^{\prime\mix}(R,G)$
is given by~$F_{j_0\dots j_m}(P,a,a')=(P,b,b')$,
with
	$$b(i_0,\dots,i_l)
	=a\(j_{i_0},\dots,j_{i_l}\)
		\Mand
	b'(i_0,\dots,i_{l-1})
	=f_{j_0\dots j_m}^*\(a(j_{i_0},\dots,j_{i_{l-1}})\)$$
for all~$0\le i_0\le\dots\le i_{l-1}\le m$.
\endroster
The {\em mixed
Whitehead space\/}~$\Wh^{\prime\mix}(R,G)$
is the geometric realisation of the mixed Whitehead category.
%
\enddefinition

\Remark\WhMixedFlatBdlRem
As in \WhFlatBdlRem\ and \WhDiffFlatBdlRem,
the complex representation~$\rho$ gives rise to flat
bundles~$\Cal V\to\Wh^{\prime\mix}(R,G)$ and~$\Cal H\to\Wh^{\prime\mix}(R,G)$.
Here, $\Cal V$ is now equipped with a mixed superconnection.
\endremark

\Example\WhMixedExp
The coefficients~$a'(\sigma,\punkt)$ constructed in \SimplPatchProp\
clearly give rise
to an object~$(P_\sigma,a_\sigma,a'_\sigma)\in\cWh^\mix_k(R;G)$
if~$\sigma\in S_k$,
by
	$$a_\sigma(i_0,\dots,i_l)
	=a(\sigma_{i_0\dots i_l})
		\Mand
	a'_\sigma(i_0,\dots,i_{l-1})
	=\sigma^*a'(\sigma,\sigma_{i_0\dots i_{l-1}})\;.$$
Mimicking the proof of~\WhFunctorProp,
we can show that this defines a functor~$\Xi_{p,h,F}^\mix\colon\Cal S
\to\cWh^{\prime\mix}(R;G)$.
\endexample

To pass to the full Whitehead category,
we define elementary extensions of differential and mixed superconnections.
Let~$P\dg j$ and~$a\dg j$ be given as in~\ElExADef.
Similarly,
we define~$a'$ as in \MixedMatrixDef\ with only one non-zero component
	$$a^{\prime[j]}_0=a^{\prime[j]}(\emptyset)
	=\Bigl(\smallmatrix0&1\\&0\endsmallmatrix\Bigr)\;.$$
We also use the symbol~``$\oplus$'' as in \ElExDef.

\Definition\DiffElExDef
The {\em elementary expansions\/} of~$(P,a')\in\cWh_k^{\prime\diff}$
and~$(P,a,a')\in\cWh_k^{\prime\mix}$
are given by~$(P,a')\oplus(P\dg j,a^{\prime[j]})$
and~$(P,a,a')\oplus(P\dg j,a\dg j,a^{\prime[j]})$,
respectively.
\enddefinition

\Definition\DiffStabDef
The {\em (stabilised) differential Whitehead category\/}~$\cWh^\diff(R,G)$
has the same objects and face operators as~$\cWh^{\prime\diff}(R,G)$.
Its morphisms are iterated formal compositions
of $\cWh^{\prime\diff}(R,G)$-morphisms and elementary expansions.
The {\em (stabilised) differential Whitehead space\/}~$\Wh^\diff(R,G)$
is the geometric realisation of~$\cWh^\diff(R,G)$.

The {\em (stabilised) mixed Whitehead category\/}~$\cWh^\mix(R,G)$
and the {\em (stabilised) mixed Whitehead space\/}~$\Wh^\mix(R,G)$
are defined analogously.

As before,
we define the full subcategories~$\cWh^{\diff,h}(R,G)\subset\cWh^\diff(R,G)$
and~$\cWh^{\mix,h}(R,G)\subset\cWh^\mix(R,G)$
consisting of all fibre-wise acyclic objects.
\enddefinition

\Remark\WhForgetfulRem
The category~$\cWh^\mix(R;G)$ encompasses both~$\cWh^\simp(R;G)$
and~$\cWh^\diff(R;G)$.
In fact, there are forgetful functors
$$\alignat4
	&&\Phi_1\colon\cWh^\mix(R;G)&\to\cWh^\simp(R;G)
		&\Mtext{with}&&
	(P,a,a')&\mapsto(P,a)\;,\\
	\Land
	&&\Phi_2\colon\cWh^\mix(R;G)&\to\cWh^\diff(R;G)
		&\Mtext{with}&&
	(P,a,a')&\mapsto\(P,a'(\emptyset)\)\;.
\endalignat$$
If~$(f,g)$ is a morphism in~$\cWh^\mix(R;G)$,
then~$\Phi_1(f,g)$ and~$\Phi_2(f,g)$ are still given by~$(f,g)$
in the corresponding category.

It follows from the proof of \SimplPatchProp\
that~$\Phi_1$ has a right-inverse on objects.
We will now show that in fact,
both functors have right-inverses as functors,
which make them homotopy equivalences on the level of realisations.
This implies in particular
that during the passage from a simplicial superconnection on~$V\to B$
to a smooth one in \ComplexKapitel,
up to homotopy,
no information was lost.
It also implies that the simplicial, mixed, and smooth superconnections
of \ComplexKapitel\ are all ``classified'' by maps into the same
pre-Whitehead space~$\Wh'(R;G)$
that was used by Igusa to construct his torsion classes~$\tau_*(M/B;F)$.
\endremark

\Theorem\WhHomEqThm.
The forgetful functors
	$$\Phi_1\colon\cWh^\mix(R;G)\to\cWh^\simp(R;G)
		\Mand
	\Phi_2\colon\cWh^\mix(R;G)\to\cWh^\diff(R;G)$$
have homotopy inverses~$\Psi_1$ and~$\Psi_2$ on realisations.

\Remark\WhFlatBdlEqRem
For each complex representation~$\rho$ of~$R$,
we clearly have
$$\align
	\Cal V\
	\cong\abs{\Phi_1}^*\Cal V
	\cong\abs{\Phi_2}^*\Cal V
	&\to\Wh^{\prime\mix}(R,G)\;,\\
		\Land
	\Cal H\
	\cong\abs{\Phi_1}^*\Cal H
	\cong\abs{\Phi_2}^*\Cal H
	&\to\Wh^{\mix}(R,G)
\endalign$$
as flat bundles.
Also by construction,
the flat simplicial superconnection on~$\Cal V$
and the flat differential superconnection on~$\Cal V$
arise by forgetting parts of the mixed superconnection
on~$\Cal V\to\Wh^{\prime\mix}(R,G)$.
\endremark

\demo{Proof\/ \rm of \WhHomEqThm}
We start with the functor~$\Phi_1$,
and we work on the level of pre-Whitehead categories,
taking care that all our constructions
are compatible with stabilisation.
The construction in the proof of \SimplPatchProp\
assigns to each object~$(P,a)$ of~$\cWh^{\prime\simp}(R;G)$
an object~$(P,a,a')$ of~$\cWh^{\prime\mix}(R;G)$.
However, this construction is not natural,
so it is not clear that we obtain a functor from~$\cWh^{\prime\simp}(R;G)$
to~$\cWh^{\prime\mix}(R;G)$.
Since morphisms are unaffected by generalised face operators
and by the functor~$\Phi_1$,
it is clear that a right-inverse~$\Psi_1$ of~$\Phi_1$
must also map each $\cWh^{\prime\simp}(R;G)$-morphism~$(f,g)$
to the $\cWh^{\prime\mix}(R;G)$-morphism~$(f,g)$.

We proceed by induction over the simplicial degree.
In degree~$0$,
the functor~$\Phi_1$ restricts to an identity of categories
	$$\cWh_0^\mix(R;G)=\cWh_0^\simp(R;G)\;,$$
so~$\Psi_{1,0}=\Phi_{1,0}^{-1}$ is simply the inverse of this identity.
We assume now
that we have constructed
functors~$\Psi_{1,k'}\colon\cWh_{k'}^\simp(R;G)\to\cWh_{k'}^\mix(R;G)$
for all~$k'<k$,
such that~$\Phi_1\circ \Psi_{1,k'}$ is the identity on~$\cWh_{k'}^\simp(R;G)$.

Given any object~$(P,a)\in\cWh_k^\simp(R;G)$,
we want to construct~$\Psi_{1,k}(P,a)\in\cWh^\mix(R;G)$.
If~$(P,a)$ is degenerated,
then~$\Psi_{1,k}(P,a)$ is already determined by compatibility
with generalised face operators.
Otherwise,
we first replace~$P$ by the graded poset~$P'$ which is identical with~$P$
as a graded set,
but carries the minimal partial ordering that supports~$a$.
In other words,
we have $\alpha\prec'\beta$
iff there exists~$0\le i_0<\dots<i_l\le k$
such that~$a_{\alpha\beta}(i_0,\dots,i_l)\ne0$.
Then it follows from \WhSimpDef~(2)
that any morphism~$(f,g)\colon(P',a)\to(Q,b)$ is invertible.
Also,
assume that~$(Q,b)$ lies in the same connected component
of the nerve~$N\cWh_k^{\prime\simp}(R;G)$,
then there exists a morphism~$(Q,b)\to(P',a)$.

Assume that we can construct
an $\Aut(P',a)$-invariant object~$(P',a,a')\in\cWh^\simp_k(R;G)$,
such that
	$$F_{j_0\dots j_{k'}}\,(P',a,a')
	=\Psi_{1,k'}\,F_{j_0\dots j_{k'}}\,(P',a)$$
for all proper face operators~$F_{j_0\dots j_{k'}}$ with~$k'<k$.
Then for any morphism~$(f,g)\colon(Q,b)\to(P',a)$,
we define
	$$\Psi_{1,k}(Q,b)=(Q,b,b')
		\Mtext{with}
	b'_{\alpha\beta}(i_0,\dots,i_{l-1})
	=g_\alpha^{-1}
		\,\bar a'_{f^{-1}(\alpha),f^{-1}(\beta)}(i_0,\dots,i_{l-1})
		\,g_\beta\;.$$
By minimality of~$P'$,
the coefficients of~$b'$ are $Q$-triangular.
Also by minimality and because all morphisms are bijective,
any two morphisms~$(f,g)$, $(f',g')\colon(Q,b)\to(P',a)$
are related by an automorphism of~$(P',a)$.
Since~$(P',a,a')$ is~$\Aut(P',a)$-invariant,
it does not matter which morphism~$(f,g)\colon(Q,b)\to(P',a)$
we choose to define~$b'$.

We first do all this only for objects~$(P',a)$
that are not elementary expansions of other objects.
Compatibility with stabilisation then forces us to extend~$\Psi_{1,k}$
in a prescribed way to all remaining objects in~$\cWh_k^\simp(R,G)$.
It is now easy to check that our definition of~$\Psi_{1,k}$
is compatible with face operators, morphisms and elementary expansions.

In order to construct~$a'$,
we proceed as in the proof of \SimplPatchProp,
keeping in mind that~$a'$ should be invariant under
all automorphisms of~$a$.
Assume that~$0\le i_0<\dots i_{l-1}\le k$ with~$i_{l-1}>l-1$,
then~$a'(i_0,\dots,i_{l-1})$ is already completely determined
by~$\Psi_{1,k'}\,F_{i_0,\dots,i_{l-1},i_{l-1}+1,\dots,k}(P',a)$.
Since
	$$\Aut(P',a)
	\subset\Aut\(F_{i_0,\dots,i_{l-1},i_{l-1}+1,\dots,k}(P',a)\)\;,$$
the coefficient~$a'(i_0,\dots,i_{l-1})$ is $\Aut(P',a)$-invariant.
Let us assume
that we already inductively constructed all~$a'(0,\dots,l'-1)$
with~$l<l'\le k+1$.
Then the restriction of~$a'(0,\dots,l-1)$ to~$\partial\Delta_{k\setminus l-1}$
is determined by induction and by~\GCSSCADef\ (for~$l\ge 1$)
or~\SimplPatchOneFormel\ (for~$l=0$),
and is also $\Aut(P',a)$-invariant.
Because~$\Aut(P',a)$ acts linearly on~$T_P(R)$,
the fixpoint set is an $\R$-subalgebra of~$T_P(R)$.
We can thus find an invariant continuation
of~$a'(0,\dots,l-1)|_{\partial\Delta_{k\setminus l-1}}$
to~$\Delta_{k\setminus l-1}$.
This completes the construction of~$(P',a,a')\in\cWh^\mix_k(R;G)$.

This completes the construction of~$\Psi_1$.
It remains to show that~$\Psi_1\circ \Phi_1$ is simplicially homotopic
to the identity on~$\cWh^\mix_k(R;G)$.
Such a homotopy can be constructed by the same procedure as above,
see e.g.\ the proofs of \RelativityProp\ and \WhFunctorProp.
Details are left to the reader.

Now, consider the functor~$\Phi_2\colon\cWh^\mix(R;G)\to\cWh^\simp(R;G)$.
The general arguments are the same as in the construction of~$\Psi_1$,
so it is enough to give an inductive $\Aut(P,a')$-invariant
construction of~$(P,\bar a,\bar a')=\Psi_1(P,a')$,
where~$(P,a')\in\cWh_k^\diff(R,G)$
is neither degenerated nor an elementary expansion,
and~$P$ carries the minimal ordering that supports~$a'$.

Note that~$\bar a'$ determines~$\bar a$ by \MixedMatrixDef~(1),
so we only construct~$\bar a'$.
As before,
all coefficients~$\bar a'(i_0,\dots,i_{l-1})$
with~$0\le i_0<\dots<i_{l-1}<k$ and~$i_{l-1}>l-1$
are inductively fixed by compatibility with face operators,
and of course,
they are $\Aut(P,a')$-invariant.
We will construct
the coefficients~$\bar a'(0,\dots,l-1)$ by induction over~$l$.
For~$l=0$, we clearly have~$\bar a'(\emptyset)=a'$.
For~$l>0$,
we already know the restriction of~$\bar a'(0,\dots,l-1)$
to the set
	$$\bigcup_{j=l}^k f_{l-1,\dots,\widehat j,\dots,k}\(\Delta_{k-l}\)
	=\partial\Delta_{k\setminus l-1}
		\setminus\inn\Delta_{k\setminus l}\;.$$
Note that this set is a deformation retract of~$\Delta_{k\setminus l-1}$.
By \WhMixedDef~(1),
we need an $\Aut(P,a')$-invariant solution~$\bar a'(0,\dots,l-1)$
of the system of linear partial first-order differential equations
	$$(-1)^l\,d\(\bar a'_{\alpha\beta}(0,\dots,l-1)\)
	=-(\del'\bar a'+a'*'a')_{\alpha\beta}(0,\dots,l-1)\;,
		\Formel\WhPsiZwoDGL$$
such that~$\bar a'(0,\dots,l-1)$
has degree~$1-l$ in~$\Omega^*\(\Delta_{k\setminus l-1};T_P(R)\)$
and restricts to the given initial data on~$\partial\Delta_{k\setminus l-1}
\setminus\inn\Delta_{k\setminus l}$.
By triangularity of~$\bar a'$,
this system can be solved inductively in the distance of~$\alpha$ and~$\beta$
in the partial ordering of~$P$.
The integrability condition is inductively given by
	$$0
	=\(d'(\del'\bar a'+\bar a'*'\bar a')\)_{\alpha\beta}(0,\dots,l-1)\;,$$
which follows from the equations
	$$0=\(d'\bar a'+\delta'\bar a'+\bar a'*'\bar a'\)
		_{\gamma\delta}(i_0,\dots,i_m)$$
for those indices~$\gamma$, $\delta$
and those sequences~$0\le i_0<\dots<i_m\le k$
where all coefficients are already known,
because
$$\align
	\kern4em&\kern-4em
	\(d'(\del'\bar a'+\bar a'*'\bar a')\)_{\alpha\beta}(0,\dots,l-1)\\
	&=\(-\del'd'\bar a'+(d'\bar a')*'\bar a'
		-\bar a'*'d'\bar a'\)_{\alpha\beta}(0,\dots,l-1)\\
	&=\(\del'(\del'\bar a'+\bar a'*'\bar a')
		-(\del'\bar a'+\bar a'*'\bar a')*'\bar a'
		+\bar a'*'(\del'\bar a'+\bar a'*'\bar a')
			\)_{\alpha\beta}(0,\dots,l-1)\\
	&=0\;.
\endalign$$

In order to get initial data for a Cauchy problem,
we have to specify all coefficients of~$\bar a'(0,\dots,l-1)$
along~$\partial\Delta_{k\setminus l-1}\setminus\inn\Delta_{k\setminus l}$,
including those that involve directions normal to this set.
It is clear that this is possible.
Taking for~$\bar a'(0,\dots,l-1)$ any $\Aut(P,a')$-invariant
solution of~\WhPsiZwoDGL,
we have finished the inductive construction
of~$(P,\bar a,\bar a')=\Psi_2(P,a')$.

This shows that~$\Phi_2$ has a right-inverse~$\Psi_2$.
By the same argument as above,
the simplicial functor~$\Psi_2\circ\Phi_2$ is simplicially homotopic
to the identity on~$\cWh^\mix(R;G)$.
\qed\enddemo

\Remark\IgusaHomEqRem
In~\ITwisted,
Igusa constructs a simplicial functor~$\cWh^\diff(R,G)\to\cWh^\simp(R,G)$
using Chen's iterated integrals.
It seems that on realisations,
Igusa's functor is homotopic to~$\Phi_1\circ\Psi_2$.
\endremark

\enddocument